\documentclass[12pt]{amsart}
\usepackage{amsfonts,amsmath,amsthm,amscd,amssymb,latexsym,cite,verbatim,texdraw,floatflt,
caption2}
\RequirePackage{hyperref}

\usepackage[a4paper]{geometry}

\theoremstyle
{plain}

\date{\today}

\address{I. Protasov, Department of Computer Science and Cybernetics, Kyiv University, Volodymyrska 64, 01033, Kyiv, Ukraine}
\email{i.v.protasov@gmail.com; }
\address{K. Protasova, Department of Computer Science and Cybernetics, Kyiv University, Volodymyrska 64, 01033, Kyiv, Ukraine}
\email{ksuha@freenet.com.ua}

\begin{document}

\title{Free coarse groups}

\author{Igor Protasov and Ksenia Protasova}

\maketitle
\vskip 5pt

{\bf Abstract}
A coarse group is a group endowed  with a coarse structure so that the group multiplication and inversion are coarse  mappings. Let $(X, \mathcal{E})$ be a coarse space and let  $\mathfrak{M}$  be a variety of groups different from the variety of singletons. We prove that there is a  coarse group $F_{\mathfrak{M}} (X, \mathcal{E})\in \mathfrak{M}$ such  that $(X, \mathcal{E}) $ is a subspace of
$F_{\mathfrak{M}} (X, \mathcal{E})$, $X$
 generates $F_{\mathfrak{M}} (X, \mathcal{E})$
  and every coarse mapping $(X, \mathcal{E}) \longrightarrow (G, \mathcal{E}^{\prime}) $
  where $G\in\mathfrak{M}$,  $(G, \mathcal{E}^{\prime}) $
  is a coarse group,
  can be extended to coarse homomorphism
 $F_{\mathfrak{M}} (X, \mathcal{E})\longrightarrow (G, \mathcal{E}^{\prime}) $.
If $\mathfrak{M}$  is the variety of all groups, the groups $F_{\mathfrak{M}} (X, \mathcal{E})$
are  asymptotic counterparts of Markov free  topological  groups over Tikhonov spaces.

\vskip 10pt

\vskip 5pt


\vskip 5pt

In \cite{b4}  A. A. Markov proved that, for every Tikhonov space
 $(X, \mathcal{T})$ there exists a group topology $\mathcal{T}^{\prime}$  on the free group
 $F(X)$ in the alphabet $X$ such that $(X, \mathcal{T})$  is a closed subset of $(F(X), \mathcal{T}^{\prime})$
 and every continuous mapping from $(X, \mathcal{T})$ to a topological group $G$ can be
 extended to continuous  homomorphism   $(F(X), \mathcal{T}^{\prime})\longrightarrow G$ .
In particular, every Tikhonov space can be embedded  as a closed  subset into some topological group.

Our purpose is to construct the natural counterparts of Markov free topological groups in category of coarse  groups and coarse homomorphisms. A coarse group is a group  endowed  with a coarse structure  in  such a way that the group multiplication and inversion are coarse mappings. All necessary facts about coarse  spaces and  coarse groups are in  sections 1, 2  and the construction of free coarse groups in section 3.

\section{Coarse structures}

Following \cite{b10}, we say that a family $\mathcal{E}$ of subsets of $X\times X$ is a {\it coarse structure} on a set $X$ if

\begin{itemize}
\item{} each $\varepsilon \in \mathcal{E}$ contains the diagonal $\vartriangle _{X}$, $\vartriangle _{X}= \{(x,x) : x \in X\}$ ; \vskip 5pt

\item{}  if $\varepsilon, \delta\in\mathcal{E}$ then $\varepsilon \circ\delta\in\mathcal{E}$  and $\varepsilon^{-1}\in\mathcal{E}$ where $\varepsilon \circ\delta = \{(x, y): \exists z ((x,z)\in\varepsilon, (z,y)\in\delta)\}, $  $ \ \varepsilon^{-1}= \{(y,x): (x,y)\in\varepsilon\}$;

\item{}  if $\varepsilon\in\mathcal{E}$ and $\bigtriangleup_{X}\subseteq \varepsilon^{\prime}\subseteq\varepsilon$ then $\varepsilon^{\prime}\in\mathcal{E}$.

\end{itemize}
\vskip 5pt

Each $\varepsilon\in\mathcal{E}$ is called an {\it entourage} of the diagonal.
A subset $\mathcal{E}^{\prime}\subseteq\mathcal{E}$ is called a {\it base} for $\mathcal{E}$ if, for
every $\varepsilon\in\mathcal{E}$ there exists $\varepsilon^{\prime}\in\mathcal{E}^{\prime}$ such that $\varepsilon\subseteq\varepsilon^{\prime}$.

The pair $(X, \mathcal{E})$ is called a {\it coarse space}. For
$x\in X$  and $\varepsilon\in\mathcal{E}$, we denote
 $B(x, \varepsilon)= \{y\in X: (x,y)\in\varepsilon \}$
and say that
 $B(x,\varepsilon)$ is a {\it ball of radius  $\varepsilon$ around $x$.}
 We note that a coarse space can be considered as
 an asymptotic counterpart of a uniform topological space and could
  be defined in terms of balls, see
 \cite{b7}, \cite{b9}.
 In this case a coarse space is called a {\it ballean}.

A coarse space $(X,\mathcal{E})$ is called {\it connected} if, for any $x, y \in X$,
 there exists  $\varepsilon\in\mathcal{E}$ such that $y\in B(x,\varepsilon)$.
 A subset $Y$ of $X$ is called {\it bounded} if there exist $x\in X$ and $\varepsilon\in\mathcal{E}$
   such that $Y\subseteq B(x, \varepsilon)$.  The coarse
   structure
    $\mathcal{E}=\{\varepsilon\in X\times X: \bigtriangleup_{X}\subseteq\varepsilon\}$ is
     the unique coarse structure such that  $(X,\mathcal{E})$ is connected and bounded.

In what follows, all coarse spaces under consideration are supposed to be {\bf connected}.

Give  a coarse space $(X, \mathcal{E})$, each subset $Y \subseteq X$ has the natural coarse structure
$\mathcal{E}|_{Y}= \{\varepsilon\cap(Y\times Y): \varepsilon\in\mathcal{E} \}$, $(Y, \mathcal{E}|_{Y})$
 is called a {\it subspace} of $(X, \mathcal{E})$.
A subset $Y$ of $X$  is called {\it large} (or {\it  coarsely dense})  if there exists $\varepsilon\in \mathcal{E}$  such that $X= B(Y, \varepsilon)$ where $B(Y, \varepsilon)=\cup_{y\in Y} B(Y, \varepsilon)$.

Let  $(X, \mathcal{E})$, $(X^{\prime}, \mathcal{E}^{\prime})$ be coarse spaces. A mapping
$f: X\longrightarrow X^{\prime}$
is called {\it coarse}  (or {\it  bornologous}  in terminology of \cite{b8}) if, for every
$\varepsilon\in\mathcal{E}$
there exists $\varepsilon^{\prime}\in\mathcal{E}$ such that, for every $x\in X$,  we have
$f(B(x,\varepsilon))\subseteq (B(f(x),\varepsilon^{\prime}))$.
If $f$ is surjective and coarse then $(X^{\prime}, \mathcal{E}^{\prime})$
 is called a {\it coarse image} of $(X, \mathcal{E})$.
If $f$ is a  bijection such that $f$  and $f^{-1}$ are coarse mappings   then $f$ is called an
 {\it asymorphism}.
The coarse spaces
$(X, \mathcal{E})$, $(X^{\prime}, \mathcal{E}^{\prime})$
 are called {\it coarsely equivalent} if there exist large subsets
 $Y\subseteq X$, $Y^{\prime}\subseteq X$ such that
 $(Y, \mathcal{E}|_{Y})$
 and $(Y^{\prime}, \mathcal{E}^{\prime}|_{Y^{\prime}})$
  are asymorphic.

To conclude the coarse vocabulary, we take a family
$\{(X_{\alpha}, \mathcal{E}_{\alpha}) :  \alpha< \kappa\}$
  of coarse  spaces and define the
  {\it product}
 $P_{\alpha< \kappa}(X_{\alpha}, \mathcal{E}_{\alpha})$
  as the set $P_{\alpha< \kappa} X_{\alpha}$
 endowed with the coarse  structure with the base
 $P_{\alpha< \kappa} \mathcal{E}_{\alpha}$.
 If $\varepsilon_{\alpha}\in\mathcal{E}_{\alpha}$, $\alpha<\kappa$
  and  $x,y\in P_{\alpha<\kappa}X_{\alpha}$,
  $x=(x_{\alpha})_{\alpha<\kappa}$, $y=(y_{\alpha})_{\alpha<\kappa}$
  then $(x,y)\in (\varepsilon_{\alpha})_{\alpha<\kappa}$
  if and only if $(x_{\alpha}, y_{\alpha})\in\varepsilon_{\alpha}$
   for every $\alpha<\kappa$.

\section{Coarse groups }

Let $G$ be a group with  the identity $e$.
For a cardinal $\kappa$, $[G]^{<\kappa}$ denotes the set $\{Y\subseteq G: |Y|<\kappa\}$.

A family $\mathcal{I}$  of subsets of $G$ is called a {\it group ideal} if $\mathcal{I}$  is closed and formations of subsets and finite unions,
 $[G]^{<\omega}\subseteq\mathcal{I}$
  and $AB^{-1}\in\mathcal{I}$
   for all $A,B\in \mathcal{I}$.

   A group ideal $\mathcal{I}$ is called {\it invariant} if
$\cup _{g\in G} \ \ g^{-1}A g\in\mathcal{I}$
 for each $A\in \mathcal{I}$.
 For example, $[G]^{<\kappa}$
  is a group ideal for any infinite cardinal $\kappa$.
If $\kappa>|G|$ we get the ideal $\mathcal{P}_{G}$ of all subsets of $G$.
We note also that $[G]^{<\omega}$ is invariant if and only if the set
$\{x^{-1}gx: x\in G\}$
is finite for each $g\in G$. By \cite{b6},  for every countable group $G$, there are $2^{2^{\omega}}$ distinct group ideals on $G$

Let $X$ be a $G$-space with the action  $G\times  X\longrightarrow  X$, $(g,x)\longmapsto  gx$.
We assume that $G$ acts on $X$ transitively,  take a group ideal $\mathcal{I}$ on $G$  and consider the coarse structure $\mathcal{E}(G, \mathcal{I}, X)$  on $X$  with  the base $\{\varepsilon_{A}: A\in\mathcal{I}, \  e\in A\}$,  $ \ \varepsilon_{A}=\{(x, gx): x\in X, g\in A\} $.
Then $(x,y)\in \varepsilon_{A}$  if and only if $yx ^{-1} \in A$ so
$B(x,\varepsilon)=Ax$,   $ \ Ax=\{gx: g\in A\}$.

By \cite[Theorem 1]{b5}, for every coarse structure  $\mathcal{E}$ on $X$,
 there exist a group $G$ of permutations of $X$ and a group ideal $\mathcal{I}$ on $G$  such that
 $\mathcal{E}= \mathcal{E} (G, \mathcal{I}, X)$.

Now let $X=G$  and $G$  acts on $X$  by the left shifts.
We denote
$\mathcal{E}_{\mathcal{I}}= \mathcal{E} (G, \mathcal{I}, G)$.
Thus, every group ideal $\mathcal{I}$ on $G$ turns $G$ into the coarse space
$(G, \mathcal{E}_{\mathcal{I}})$.
We note that a subset $A$ of $G$ is bounded in $(G, \mathcal{E}_{\mathcal{I}})$ if and only if  $A\in\mathcal{I}$.

For finitely generated groups, right coarse groups $(G, \mathcal{E}_{[G]<\omega})$ in metric form take a great part of
{\it Geometrical Group Theory}, see \cite[Chapter 4]{b2}.

A group $G$ endowed with a coarse structure $\mathcal{E}$  is called {\it left (right) coarse group} if, for every $\varepsilon\in \mathcal{E}$,  there exists $\varepsilon^{\prime}\in \mathcal{E}$
 such that $g B(x,\varepsilon)\subseteq B(gx,\varepsilon^{\prime})$ $((B(x,\varepsilon)g \subseteq B(xg,\varepsilon^{\prime}))$
  for all $x, g\in G$.

A group $G$  endowed with a coarse structure $\mathcal{E}$  is called a
{\it coarse group} if the group multiplication $(G, \mathcal{E})\times (G, \mathcal{E})\longrightarrow (G, \mathcal{E}),$ $ \ (x,y)\longmapsto xy$  and the   inversion
$(G, \mathcal{E})\longrightarrow (G, \mathcal{E}), $ $ \ x\longmapsto x^{-1}$
  are coarse mappings. In this case, $\mathcal{E}$ is called a group coarse structure.

For proofs of the following two statements see \cite{b8} or \cite[Section 6]{b9}.
\vskip 6pt

{\bf Proposition 1.  } {\it
A group $G$ endowed with a coarse structure  $\mathcal{E}$ is a right coarse group if and  only if there exists a group ideal $\mathcal{I}$ on $G$ such that $\mathcal{E}=\mathcal{E}_{\mathcal{I}}$. }
\vskip 6pt

{\bf Proposition 2.  } {\it
For a group $G$ endowed with a coarse structure  $\mathcal{E}$, the following conditions are equivalent:

\vskip 5pt

(i) 	$(G,\mathcal{E})$  is a coarse group;
\vskip 5pt

(ii)   $(G,\mathcal{E})$ is  left and right coarse  group;
\vskip 5pt

(iii)  there exists an invariant group ideal $\mathcal{I}$ on $G$  such that $\mathcal{E}=\mathcal{E}_{\mathcal{I}}$.
}
\vskip 7pt

{\bf Proposition 3.  } {\it
Every  group coarse structure $\mathcal{E}$ on a subgroup $H$ of an Abelian group $G$ can be extended to a group coarse structure $\mathcal{E}^{\prime}$ on $G$.\vskip 5pt

Proof.}
We take a group ideal $\mathcal{I}$ on $G$ such that $\mathcal{E}=\mathcal{E}_{\mathcal{I}}$, denote by
$\mathcal{I}^{\prime}$   the group ideal on $G$ with the base $A+B$,
$A\in [G] ^{<\omega}$, $B\in \mathcal{I}$  and  put  $\mathcal{E}^{\prime}=\mathcal{E}_{\mathcal{I}^{\prime}}$.
$ \  \  \  \Box$

\vskip 7pt

{\bf Example 1.}  We construct a group  $G$ with a normal Abelian
subgroup $H$
 of index $|G: H|=2$  such that some group
 coarse structure $\mathcal{E}$  on $H$ can not be extended to a  right
  group coarse structure on $G$. Let $H=\otimes_{n\in\mathbb{Z}} C_{n}$,
  $C_{n}\simeq\mathbb{Z}_{2}$.
  Every element $a\in H$  can be written as
  $a=(a_{n})_{n\in \mathbb{Z}}$   with $a_{n}\in C_{n}$
   and  $a_{n} =0 $ for all but finitely many $n$. We define an automorphism $\varphi$  of
   order 2 of  $H$ by  $\varphi(a_{n})_{n\in \mathbb{Z}}= (c_{n})_{n\in \mathbb{Z}}$, $c_{n}= a_{-n}$
   for each $n\in \mathbb{Z}$.
We put  $<\varphi>= \{\varphi, id\}$  and consider the
semidirect product
 $G= H \  \  \leftthreetimes <\varphi>$.
If $(h_{1}, \varphi_{1}), (h_{2}, \varphi_{2})\in G$
 then $(h_{1}, \varphi_{1}),  \  (h_{2}, \varphi_{2})= (h_{1}, \  \varphi(h_{2}), \  \varphi_{1} \varphi_{2}) $.
For each $m\in\mathbb{Z}$
we set  $H_{m}=\otimes_{n\geq m} H_{n}$.
Then the family
$\{H_{m}: m\in \mathbb{Z}\}$
is a  base  for some group ideal $\mathcal{I}$ on $G$.
We put  $\mathcal{E}=\mathcal{E}_{\mathcal{I}}$  and  take  an arbitrary invariant
group ideal $\mathcal{J}$ on $G$  such that $\mathcal{I} \subset \mathcal{J}$.
Since $\varphi H_{0}\varphi\cup  H_{0}=H$,
we see that $H\in \mathcal{J}$.
It follows that the coarse structure $\mathcal{E}_{\mathcal{J}}|_{H}$
 is bounded so $\mathcal{E}_{\mathcal{J}}|_{H}\neq\mathcal{E}$.
\vskip 7pt

 {\bf Example 2.  }
Let $G$   be an infinite group with only
two classes of conjugated elements, see \cite{b3}.
Then there is only one group coarse structure $\mathcal{E}$ on $G$,
namely $\mathcal{E}= \mathcal{E} _{\mathcal{P}(G)}$.

\section{Free coarse groups}

A class $\mathfrak{M}$ of groups is called a {\it variety}  if $\mathfrak{M}$ is
closed under formation of subgroups, homomorphic images and products.
 We assume that $\mathfrak{M}$ is non-trivial (i.e.
 there exists  $G\in \mathfrak{M}$  such that  $|G|>1$)  and recall that
 the {\it free group}  $F_{\mathfrak{M}} (X)$  is defined by the following
 conditions: $F_{\mathfrak{M}} (X)\in\mathfrak{M}$, $X\subset F_{\mathfrak{M}} (X)$,
 $X$   generates  $F_{\mathfrak{M}} (X)$
 and every mapping $X \longrightarrow G$,  $G\in \mathfrak{M}$ can be extended to  homomorphism
 $F_{\mathfrak{M}} (X)\longrightarrow G$.

Let  $(X, \mathcal{E})$ be a coarse space. We assume that
$(F_{\mathfrak{M}} (X), \mathcal{E}^{\prime})$ is a coarse group
  such that $(X, \mathcal{E})$
   is a subspace of
   $(F_{\mathfrak{M}} (X), \mathcal{E}^{\prime})$
   and every coarse mapping
   $(X, \mathcal{E})\longrightarrow (G, \mathcal{E}^{\prime\prime})$,
   $G\in \mathfrak{M}$, $(G, \mathcal{E^{\prime\prime}})$
   is a coarse group, can be extended to  coarse homomorphism
   $(F_{\mathfrak{M}} (X), \mathcal{E}^{\prime})\longrightarrow (G, \mathcal{E^{\prime\prime}})$.
We observe that this $\mathcal{E}^{\prime}$  is unique, denote
$F_{\mathfrak{M}} (X, \mathcal{E})=(F_{\mathfrak{M}} (X), \mathcal{E}^{\prime})$
and say that  $F_{\mathfrak{M}} (X, \mathcal{E})$
is a {\it free coarse group} over $(X, \mathcal{E})$ in the variety $\mathfrak{M}$.

Our goal is to prove the existence of $F_{\mathfrak{M}} (X, \mathcal{E})$
 for every coarse space $(X, \mathcal{E})$  and every non-trivial variety $\mathfrak{M}$.
\vskip 7pt

 {\bf Lemma 1.  } {\it Let  $(X, \mathcal{E})$ be a coarse space. If there is a group coarse structure
  $\mathcal{E}^{\prime}$ on $F_{\mathfrak{M}} (X)$
    such that $\mathcal{E}^{\prime}|_{X}= \mathcal{E}$
     then there exists  $F_{\mathfrak{M}} (X, \mathcal{E})$ .

\vskip 6pt

Proof.} We denote $\mathfrak{F}=\{\mathcal{T}: \mathcal{T}$  is a group coarse structure  on
$F_{\mathfrak{M}} (X)$
 such that
  $\mathcal{T}|_{X}=\mathcal{E}\}$.

By the assumption, $\mathcal{E}^{\prime}\in\mathfrak{F}$.
We set  $\mathcal{T}^{\prime}= \cap\mathfrak{F}$
  and note that $\mathcal{T}^{\prime}$
  is a group coarse structure and
  $\mathcal{T}^{\prime}|_{X}=\mathcal{E}$.
Let $G\in\mathfrak{M}$,  $(G, \mathcal{E}^{\prime\prime})$
be a coarse group,
$f: (X, \mathcal{E})\longrightarrow (G, \mathcal{E}^{\prime\prime})$
 be a coarse mapping. We extend  $f$ to homomorphism
 $f: F_{\mathfrak{M}} (X)\longrightarrow G$.
Then the coarse structure on $ F_{\mathfrak{M}} (X)$
with the base
$\{f^{-1}(\varepsilon^{\prime\prime}): \varepsilon^{\prime\prime}\in\mathcal{E}^{\prime\prime} \}$
is in $\mathfrak{F}$.
It follows that  the homomorphism
$f: (F_{\mathfrak{M}} (X), \mathcal{T}^{\prime}) \longrightarrow (G, \mathcal{E}^{\prime\prime})$
 is coarse. Hence,
 $ (F_{\mathfrak{M}} (X), \mathcal{T}^{\prime}) = \mathfrak{F}_{\mathfrak{M}}
  (X, \mathcal{E})$.$ \  \  \  \Box$

\vskip 7pt

{\bf Lemma 2.} {\it
For every coarse space  $(X, \mathcal{E})$  and every non-trivial variety $\mathfrak{M}$ of groups, there exists a group coarse structure $\mathcal{E}^{\prime}$  on $F_{\mathfrak{M}} (X)$  such  that $\mathcal{E}^{\prime}|_{X} = \mathcal{E}$.
\vskip 6pt

Proof.}
For some prime number  $p$, $\mathfrak{M}$ contains the variety
$\mathcal{A}_{p}$ of all Abelian group of period $p$.
We prove the theorem for $\mathcal{A}_{p}$  and then for $\mathfrak{M}$.

We take the free group $A(X) $  over $X$  in $A_{p}$.
Every non-zero element $a\in  A(X)$  has the unique (up to permutations of items) representation
$$(\ast)  \  m_{1}x_{1} + m_{2}x_{2} + \ldots + m_{k}x_{k},  \  \  x_{i}\in X,  \ m_{i}\in \mathbb{Z}_{p}\setminus \{0\}, \  \   i\in\{1, \ldots , k\}.$$

For every $\varepsilon\in \mathcal{E}$, $\varepsilon= \varepsilon^{-1}$
 we denote
$Y_{\varepsilon}=\{x-y:  x, y \in X,   (x,y)\in\varepsilon\}$
and by $Y_{n,\varepsilon}$ the sum on $n$ copies of   $Y_{\varepsilon}$.
We take $z\in X$ and consider the ideal $\mathcal{I}$  on $A(X)$  with the base
$$Y_{n,\varepsilon} + \{0, z, 2z, \ldots , (p-1)z\},  n<\omega.$$

We  note that $Y_{n,\varepsilon} - Y_{n^{\prime},\varepsilon^{\prime}} \subseteq   Y_{n+n^{\prime},\varepsilon\circ\varepsilon^{\prime}} $.
It follows that $B-C\in\mathcal{I}$
 for all $B, C\in\mathcal{I}$.
To show that $[F_{\mathfrak{M}}(X)]^{<\omega}\in\mathcal{I}$,
we take $x\in X$  and find $\varepsilon\in\mathcal{E}$ such that
$(x,z)\in \varepsilon$.
Then $x-z\in Y_{\varepsilon}$
and  $x\in Y_{\varepsilon} + z$.
Hence, $\mathcal{I}$   is  a group ideal.
We put $\mathcal{E}^{\prime}=\mathcal{E}_{\mathcal{I}}$ and show that   $\mathcal{E}^{\prime}|_{X}=\mathcal{E}$.

If $\varepsilon\in\mathcal{E}$, $\varepsilon=\varepsilon^{-1}$ and $(x,y)\in\varepsilon$
 then $x-y\in Y_{\varepsilon}$ so $\mathcal{E}\subseteq\mathcal{E}^{\prime}$. To prove the inverse inclusion, we take
$Y_{n,\varepsilon}+\{0, z, \ldots , (p-1)z\}$,
 assume that $x-y\in Y_{n,\varepsilon}+ \{0, z, \ldots , (p-1)z\}$
   and consider two cases.
\vskip 7pt

{\it Case:} $x-y\in Y_{n,\varepsilon}+ iz,$  $ \ i\neq 0$.
We denote by $H$  the subgroup of all  $a\in A(X)$  such that
$m_{1} + \ldots + m_{k} =0 (mod \ p)$ in the canonical representation $( \ast)$.
Then $x-y\in H$,  $Y_{n,\varepsilon}\subseteq H $ but $iz\notin H$  so this case is impossible.

\vskip 7pt

{\it Case:} $x-y\in Y_{n,\varepsilon}.$
We show that $(x,y)\in\varepsilon^{n}$.
We write $x-y$  as $(x_{1}- y_{1})+ \ldots + (x_{n}- y_{n})$, $x_{i},y_{i}\in Y_{\varepsilon}$
   so   $(x_{i},y_{i})\in\varepsilon$.
Assume that there exists $k\in\{1, \ldots, n-1\}$
such that
$$\{x_{1},y_{1}, \ldots , x_{k},y_{k}\}\cap \{x_{k+1},y_{k+1}, \ldots , x_{n},y_{n}\}=\emptyset . $$

Then either
$(x_{1}- y_{1})+ \ldots + (x_{k}- y_{k})=0$, or
$(x_{k+1}- y_{k+1})+ \ldots + (x_{n}- y_{n})=0$.
Otherwise, $x-y$  in the representation  $( \ast)$ has  more
 then two items.  It follows that there is a representation
 $$x-y=(x_{1}^{\prime}- y_{1}^{\prime})+ \ldots + (x_{k}^{\prime}- y_{k}^{\prime}), \ \ x_{i}^{\prime}, y_{i}^{\prime}\in Y_{\varepsilon}, \ \ i\in\{1, \ldots, k\}, \ \   k\leq n $$
  such that
  $\{x_{i+1}^{\prime}, y_{i+1}^{\prime}\}\cap \{x_{1}^{\prime}, y_{1}^{\prime}, \ldots , x_{i}^{\prime}, y_{i}^{\prime}\}\neq\emptyset$
   for each  $i\in\{1,\ldots, k-1\}$.
If $(x^{\prime}, y^{\prime})\in\varepsilon^{i}$ for all $x^{\prime}, y^{\prime}\in \{x_{1}^{\prime}, y_{1}^{\prime}, \ldots, x_{i}^{\prime}, y_{i}^{\prime} \}$
 then $(x^{\prime}, y^{\prime})\in\varepsilon^{i+1}$
 for all
 $x^{\prime}, y^{\prime}\in \{x_{1}^{\prime}, y_{1}^{\prime},
  \ldots, x_{i+1}^{\prime}, y_{i+1}^{\prime} \}$.
After $k$ steps, we get $x-y \in  \varepsilon^{k}$  so  $x-y\in \varepsilon^{n}$.

To  conclude the proof, we extend the mapping $id: X\longrightarrow  X$ to  homomorphism
$f: F_{\mathfrak{M}}(X)\longrightarrow A(X)$.
Then  $\{f^{-1}(Y): Y\in\mathcal{I}\}$
 is a base for some invariant group ideal $\mathcal{J}$ on $F_{\mathfrak{M}}(X)$.
Then  $(F_{\mathfrak{M}}(X), \mathcal{E}_{\mathcal{J}})$
  is a coarse group. Since $f|_{X}= id$,
  we  have $\mathcal{E}_{\mathcal{J}}|_{X}= \mathcal{E}_{\mathcal{J}}|_{X}= \mathcal{E}$.
  $ \  \  \  \Box$

\vskip 10pt

{\bf Theorem.} {\it
For every coarse space $(X, \mathcal{E})$  and every non-trivial  variety $\mathfrak{M}$ of groups,  there exists the free coarse group $F_{\mathfrak{M}}(X, \mathcal{E})$.

\vskip 6pt

Proof.} Apply Lemma 2 and Lemma 1.$ \  \  \  \Box$

\vskip 10pt

{\bf Remark  1.}
To describe the coarse structure $\mathcal{E}^{\ast}$  of  $ F_{\mathfrak{M}} (X, \mathcal{E})$ explicitly,
 for every  $\varepsilon\in \mathcal{E}$, we put  $\mathcal{D}_{\varepsilon}=\{xy^{-1}: x,y\in X, \  \  (x,y)\in\varepsilon\}$,  take $z\in X$  and denote by $P _{n,\varepsilon}$  the product on $n$ copies of the set
$$\bigcup_{g\in F_{\mathfrak{M}}(X)} \ \ g^{-1} (\mathcal{D}_{n, \varepsilon} \bigcup \mathcal{D}_{n, \varepsilon} \ z) \ g .$$
Then
$\{P_{n,\varepsilon}: \varepsilon\in\mathcal{E}, n<\omega\}$
is a base for some invariant group ideal  $\mathcal{I}^{\ast}$ on  $F_{\mathfrak{M}} (X)$.
Each subset $A\in I^{\ast}$ is bounded  in $F(X, \mathcal{E})$ so
$\mathcal{E} _{\mathcal{I}^{\ast}} \subseteq \mathcal{E}^{\ast}$.  To see that
$\mathcal{E}^{\ast}  \subseteq \mathcal{E} _{\mathcal{I}^{\ast}} $,
 the reader can repeat the arguments concluding the proof of  Lemma 2.
Hence,
$\mathcal{E}^{\ast}  = \mathcal{E} _{\mathcal{I}^{\ast}} $
\vskip 10pt

{\bf Remark 2.}  Each metric space  $(X, d)$
 defines the coarse structure
 $\mathcal{E}_{d}$ on $X$ with the base
 $\{(x,y): d(x,y)<n\}$, $n<\omega$.
By \cite[Theorem 2.1.1]{b9},  a coarse   structure $\mathcal{E}$ is metrizable if  and only if $\mathcal{E}$  has a countable  base.  If $\mathcal{E}$  is  metrizable then, in view of Remark 1,  the coarse  structure of
$F _{\mathfrak{M}} (X, \mathcal{E})$   is  metrizable.
\vskip 10pt

{\bf  Remark 3.}  If the coarse spaces  $(X, \mathcal{E}), (X, \mathcal{E}^{\prime})$  are asymorphic
 then evidently $F _{\mathfrak{M}} (X, \mathcal{E})$, $F _{\mathfrak{M}} (X^{\prime}, \mathcal{E}^{\prime})$
    are asymorphic  but this is not true with coarse equivalences in place of asymorphisms.

Let $\mathfrak{M}= \mathcal{A} _{p}$ and let $X$  be an infinite set endowed
with the bounded coarse structure $\mathcal{E}$.
We take
 $X^{\prime},  |X^{\prime}|=1$  and denote   by $\mathcal{E}^{\prime}$  the unique  coarse structure on $X^{\prime}$.
Clearly,
$(X, \mathcal{E})$  and $(X^{\prime}, \mathcal{E}^{\prime})$
 are coarsely equivalent, and
  $F _{\mathfrak{M}} (X^{\prime}, \mathcal{E}^{\prime})$
  is a cyclic group of order $p$ with bounded coarse structure.
  To see that
  $F _{\mathfrak{M}} (X)$  is unbounded,
   we take the subset $Y _{n,\varepsilon}$  (see proof of Lemma 2)
   and note that the length of any element from  $Y _{n,\varepsilon}$   in representation $(\ast)$
   does not exceed $2n$,  but  $F _{\mathfrak{M}} (X)$  has  elements of any length.
\vskip 10pt

{\bf  Remark 4.}
Let $X$  be a Tikhonov space with distinguished point $x_{0}$.
In  \cite{b1}  M. I. Graev  defined a group topology on
$F(X\setminus \{x_{0}\} )$  such that $X$  is a closed subset of
$F(X\setminus \{x_{0}\} )$,  $x_{0}= e$,  and every continuous mapping $f: (X)\longrightarrow G $,
$f(x_{0})=e$,   $G$  is a topological group, can be extended to continuous homomorphism
$F(X\setminus \{x_{0}\} )\longrightarrow G$.

Let $(X, \mathcal{E})$  be a  coarse space with  distinguished point
$x_{0}$,  $Y= X\setminus\{x_{0}\}$,   $\mathcal{E}^{\prime}= \mathcal{E} \mid _{Y}$.
We take the free coarse group
$F(Y, \mathcal{E}^{\prime})$
and note that $\{e\}\cup Y$  is asymorphic to  $(X, \mathcal{E})$  via the mapping
 $h(y)=y$,  $y\in Y$  and $h(e)=x_{0}$.
Hence, it does not make sense to define the coarse counterparts of the Graev free topological groups.

\end{document}